\documentclass[12pt]{amsart}
\usepackage[arrow]{xy}

\theoremstyle{plain}
\begingroup
\newtheorem{thm}{Theorem}[section]
\newtheorem{lem}[thm]{Lemma}
\newtheorem{cor}[thm]{Corollary}

\endgroup

\theoremstyle{definition}
\newtheorem{defn}[thm]{Definition}

\theoremstyle{remark}
\newtheorem{rem}[thm]{Remark}
\newtheorem{ex}[thm]{Example}

\let\cal=\mathcal

\begin{document}

\def\headline #1{\medskip \centerline{\bf #1} \medskip}

\def\Ker {\operatorname{Ker}\nolimits}
\def\Im {\operatorname{Im}\nolimits}
\def\Image {\operatorname{Image}\nolimits}
\def\Syz {\operatorname{Syz}\nolimits}
\def\initial {\operatorname{in}\nolimits}
\def\Gin{\operatorname{Gin}\nolimits}
\def\Spec{\operatorname{Spec}\nolimits}
\def\D{\operatorname{D}\nolimits}
\def\V{\operatorname{V}\nolimits}
\def\H {\operatorname{H}\nolimits}
\def\E {\operatorname{E}\nolimits}
\def\V {\operatorname{V}\nolimits}
\def\nE {\operatorname{e}\nolimits}
\def\nV {\operatorname{v}\nolimits}
\def\C {\operatorname{\cal C}\nolimits}

\title
[Bipartite graphs whose edge algebras are complete intersections]
{Bipartite graphs whose edge algebras are complete intersections}
\author{Mordechai Katzman}
\address{
Departamento de matem\'aticas\hfill \\
Centro de Investigaci\'on y Estudios Avanzados
}
\email{katzman@math.cinvestav.mx}
%\keywords{edge algebras, complete intersections, toric varieties}

\begin{abstract}
Let $R$ be monomial sub-algebra of $k[x_1,\dots,x_N]$ generated by
square free monomials of degree two.
This paper addresses the following question: when is
$R$ a complete intersection?

For such a $k$-algebra we can associate a graph $G$ whose vertices are
$x_1,\dots,x_N$ and whose edges are
$\{ (x_i, x_j) | x_i x_j \in R \}$.
Conversely,
for any graph $G$ with vertices $\{x_1,\dots,x_N\}$
we define the {\it edge algebra associated with $G$} as the
sub-algebra of $k[x_1,\dots,x_N]$ generated by the monomials
$\left\{ x_i x_j   | (x_i,x_j) \text{ is an edge of } G \right\} .$
We denote this monomial algebra by $k[G]$.

This paper  describes all bipartite graphs whose edge algebras are
complete intersections.
\end{abstract}

\maketitle

%%%%%%%%%%%%%%%%%%%%%%%
\section {Introduction}
%%%%%%%%%%%%%%%%%%%%%%%

For any graph $G$ with vertices $\{x_1,\dots,x_N\}$
we define the {\it edge algebra associated with $G$} as the
sub-algebra of $k[x_1,\dots,x_N]$ generated by the monomials
$$\left\{ x_i x_j   | (x_i,x_j) \text{ is an edge of } G \right\} .$$
We denote this monomial algebra by $k[G]$.

There has been a recent effort to relate the algebraic properties of
$k[G]$ with the structure of $G$. For example, \cite{SVV}  and \cite{HO2}
give a criterion for the normality of
$k[G]$ and the authors of the latter recently obtained a characterization
of all
bipartite graphs whose edge algebras are Gorenstein (\cite{HO3}.)

In this paper I follow this line of inquiry and I will present
a characterization of all bipartite graphs whose edge algebras are
complete intersections (theorem \ref{Main-theorem}.)

We shall denote by $\E(G)$ the set of edges of $G$ and by
$\V(G)$ its set of vertices. The cardinality of these sets will be denoted
by $\nE(G)$ and $\nV(G)$ respectively.

We fix the following presentation of $k[G]$:
define a map $\phi : k[\E(G)] \rightarrow k[G]$ by
$\phi(x_r,x_s)=x_{r} x_{s}$ and let $K_G$ be the kernel of $\phi$.
Note that $K_G$ is a binomial prime ideal containing no monomials.
We also recall that if $G$ is connected then
$\dim(k[G])$ equals $\nV(G)-1$ if $G$ is bipartite,
$\nV(G)$ otherwise (corollary 7.3.1 in \cite{V}.)

The next section will describe a set of generators and a
Gr\"obner basis for $K_G$. We will then obtain a number of
straightforward corollaries linking the structure of $G$ with that
of $k[G]$.
We will then restrict our attention to bipartite graphs, and we will
characterize those graphs whose edge algebras are complete intersections.

%%%%%%%%%%%%%%%%%%%%%%%%%%%%%%%%%%%%%%
\section {A Gr\"obner basis for $K_G$}
%%%%%%%%%%%%%%%%%%%%%%%%%%%%%%%%%%%%%%

We first introduce some graph theoretical terminology:

Let $G$ be a graph. A {\it walk of length $l$ in $G$} is a sequence of edges
$$(v_1,v_2), (v_2,v_3), \dots, (v_{l-1}, v_l), (v_l, v_{l+1}) ;$$
this walk is {\it closed} if $v_{l+1}=v_1$; if, in addition,
$l\ge 3$ and every vertex in the walk occurs precisely twice
this closed path is a
{\it cycle of length $l$.}

A closed walk $e_1,e_2,\dots,e_l$ is {\it minimal} if no two consecutive (modulo $l$)
edges are equal. A walk $e_1,e_2,\dots,e_{2l}$ is {\it trivial} if after
a cyclic permutation of the edges we have
$e_1=e_2, e_3=e_4, \dots e_{2l-1}=e_{2l}$.

A closed walk $e_1^\prime , \dots , e_r^\prime$ is {\it contained in} a closed
walk $e_1,\dots,e_s$ if after a cyclic permutation of the edges of the walks
we have $e_1=e_1^\prime , \dots, e_r=e_r^\prime$.
All other graph theoretical terminology in this paper conforms with \cite{H}.

If we fix any monomial order in $k[\E(G)]$ then given any closed walk
of even length $w=e_1,\dots,e_{2l}$ we define
$$\psi(w)=\psi^+(w)-\psi^-(w)=\prod_{i=1}^l e_{2i-1} - \prod_{i=1}^l e_{2i}\in
k[\E(G)] $$
where $\psi^+(w)\ge \psi^-(w)$.
It is not hard to see that $\psi(w)\in K_G$ for all closed walks $w$
of even length and it turns out that these generate $K_G$
(lemma 1.1 in \cite{HO1}.)

\begin{thm}
Fix any lexicographic monomial order in $k[\E(G)]$.
Let $\cal W$ be the set of minimal closed walks in $G$ of even length and let
$\cal{G}_G=\{ \psi(w) | w\in {\cal W} \}$. Then there exists a subset
of $\cal{G}_G$ which is a Gr\"obner basis for $K_G$.
\end{thm}

\begin{proof}
It is enough to show that any binomial in $K_G$ reduces to zero
with respect to $\cal{G}_G$. Pick a counterexample $A-B\in K_G$
with $A>B$
of minimal degree having disjoint support.
Let $e_1$ be the largest variable occurring  in
$A$. If $e_1=(v_1,v_2)$ then some variable $e_2=(v_2,v_3)$ must occur in
$B$. If $v_3=v_1$ then $A-B=e_1(A/e_1-B/e_1)$ where the second factor is
a binomial in $K_G$ of smaller degree and by the minimality of the degree
of $A-B$ it reduces to zero, and we are done.

Assume now that $v_3\ne v_1$. We can now pick a variable $e_3=(v_3,v_4)$
occurring in $A/e_1$ and a variable $e_4=(v_4,v_5)$ occurring in $B/e_2$.
If $v_5=v_1$ then $\psi^+(e_1,e_2,e_3,e_4)$ divides $A$ and we are done.
We may continue in this fashion until we produce a closed walk
$w=e_1,e_2,\dots,e_{2l}$ such that $\psi^+(w)$ divides $A$.
\end{proof}

\begin{cor}
If $G$ has at most one odd cycle then
$$\cal{G}_G = \{ \psi(c) | c \text{ is a even cycle in } G \}$$
is a  Gr\"obner basis for $K_G$.
\end{cor}

\begin{proof}
It is enough to show that if $w=(e_1,\dots,e_{2l})$ is a minimal walk
in $G$ then $\psi^+(c)$ divides $\psi^+(w)$ for some even cycle $c$
contained in $w$. Pick as a counterexample
such a $w$ with minimal length.
Since $w$ is minimal, there exists some cycle $c^\prime$
contained in $w$, say $c^\prime=(e_1,e_2,\dots,e_r)$. If $r$ is odd then
$(e_{r+1},e_{r+2},\dots,e_{2l})$ is a closed walk of odd length, and,
therefore, cannot be trivial and must contain an even cycle.

We have shown that $w$ must contain an even cycle, say
$c=(e_1,e_2,\dots,e_{2s})$. If $s=l$ we are done,
otherwise let $w^\prime$ be the even cycle
$(e_{2s+1},e_{2s+2},\dots,e_{2l})$.
$\psi^+(w)$ must be divisible by $\psi^+(c)$ or by $\psi^+(w^\prime)$.
If the former occurs we are done, if the latter occurs, the minimality
of the length of $w$ implies that there exists an even
cycle $c^\prime$ in $w^\prime$ such that
$\psi^+(c^\prime)$ divides $\psi^+(w^\prime)$.
\end{proof}

\begin{cor}
Let $G$ have at most one odd cycle, and let $B_1,\dots,B_r$ be the blocks
of $G$.
\begin{enumerate}
\item $k[G]$ is a complete intersection if and only if
$k[B_i]$ is a complete intersection for all $1\le i \le r$.
\item $k[G]$ is Gorenstein if and only if
$k[B_i]$ is Gorenstein for all $1\le i \le r$.
\end{enumerate}
\end{cor}
\begin{proof}
Since $K_G$ is generated by elements involving edges in one block
we can write
$$k[G]\cong k[\E(G)]/K_G \cong
k[\E(B_1)]/K_{B_1} \otimes_k \dots \otimes_k k[\E(B_r)]/K_{B_r} $$
proving (1).

We can find a system of parameters for  $k[\E(G)]/K_G$ where each
parameter is in some $k[B_i]$. Killing these parameters gives us
a zero-dimensional $k$-algebra whose socle is the tensor product
of $r$ non-zero vector spaces. Thus the type of $k[G]$ is one if and only if
all these vector spaces are one dimensional.
\end{proof}

\begin{cor}
Let $G$ be a connected graph and let $d=\dim k[G]$.
Let $(1,h_1,h_2,\dots)$ be the $h$-vector of $k[G]$.
If $2L$ is the length of the smallest minimal even closed walk in $G$
then
$\displaystyle h_i=\binom{\nE(G)-d+i-1}{\nE(G)-d-1}$ for all $0\le i<L$ while
$\displaystyle \binom{\nE(G)-d+L-1}{\nE(G)-d-1}-h_L$
is the number of (minimal) closed walks
of length $2L$ in $G$.
\end{cor}

\begin{proof}
Let $\H(i)$ be the Hilbert function of $k[\E(G)]/K_G$ (where the degrees
of the variables are one,) and consider the short exact sequence
$$0\rightarrow K_G \rightarrow K[\E(G)] \rightarrow
k[\E(G)]/K_G \rightarrow 0 .$$
Since the minimal degree of a generator of $K_G$ is $L$ we have
$\H(i)=\binom{\nE(G)+i-1}{\nE(G)-1}$ for all $i<L$ while
$\H(L)=\binom{\nE(G)+L-1}{\nE(G)-1}-\gamma$ where $\gamma$ is the number of
closed walks of length $\gamma$ in $G$.

Now $h_i$ is the coefficient of $t^i$ in $(1-t)^d \sum_{j=0}^i \H(j)t^j$,
and for $i<L$ this is the coefficient of $t^i$ in
$\displaystyle \frac{(1-t)^d }{(1-t)^{\nE(G)}}$, i.e.,
$h_i=\displaystyle \binom{\nE(G)-d+i-1}{\nE(G)-d-1}$.
On the other hand, $h_L$ is the coefficient of $t^L$ in
$$
(1-t)^d \left( \sum_{j=0}^L \binom{\nE(G)+j-1}{\nE(G)-1} t^j -
 \gamma t^L \right)=\\
\frac{(1-t)^d }{(1-t)^{\nE(G)}} - (1-t)^d \gamma t^L
$$
and, therefore, $\displaystyle h_L=\binom{\nE(G)+L-1}{\nE(G)-1}-\gamma$.
\end{proof}

%%%%%%%%%%%%%%%%%%%%%%%%%%%%%%%%%%%%%%%%%%%%%%%%%%%%%%%%%%%%%%%%%%%%%%%%
\section{Bipartite graphs whose edge algebra is a complete intersection}
%%%%%%%%%%%%%%%%%%%%%%%%%%%%%%%%%%%%%%%%%%%%%%%%%%%%%%%%%%%%%%%%%%%%%%%%

We begin this section by producing a minimal set of generators
for $k[G]$ where $G$ is bipartite (i.e., all cycles in $G$ are even.)
We shall assume that we fixed some
unspecified monomial order in $k[\E(G)]$ so that for any closed walk $w$
of even length in $G$, $\psi(w)$ is well defined.

\begin{defn}
A bipartite graph $G$ is a {\it CI graph}
if any two cycles with no chords have at
most one edge in common.
\end{defn}

For any graph $G$ we will denote the set of cycles in $G$ with no chords
by $\C(G)$.

The following observation, also proved in \cite{S}, provides a
link between the structures of $G$ and $k[G]$.

\begin{thm}\label{CI-mingens}
If  $G$ is  a bipartite graph then
$${\cal S}=\{ \psi(c) | c\in \C(G) \}$$
is a minimal set of generators of $K_G$.
\end{thm}

\begin{proof}
We first show that $\cal S$ generates $K_G$. Pick as a counterexample
a cycle $c=e_1,\dots,e_{2l}$ of minimal length such that
$\psi(c)$ is not contained
in $< {\cal S} >$. Then $c$ must have a chord $e$ and we obtain after
a cyclic permutation of the edges of $c$
two cycles $c_1=e,e_1,\dots,e_{2r-1}$ and $c_2=e,e_{2r},\dots,e_{2l}$
in $G$ whose length is smaller than the length of $c$. By the minimality
of $c$ we have $\psi(c_1),\psi(c_2) \in  < {\cal S} >$ but since
\begin{gather*}
e_2 e_4 \dots e_{2r-2}
(e e_{2r+1} e_{2r+3} \dots e_{2l-1} - e_{2r} e_{2r+2} \dots e_{2l})-\\
e_{2r+1} e_{2r+3} \dots e_{2l-1}
(e e_2 e_4 \dots e_{2r-2} - e_1 e_3 \dots e_{2r-1})=\\
e_1 e_3 \dots e_{2r-1} e_{2r+1} \dots e_{2l-1} -
e_2 e_4 \dots e_{2r-2} e_{2r} \dots e_{2l}
\end{gather*}
$\psi(c)$ is in the ideal generated by $\psi(c_1)$ and $\psi(c_2)$,
a contradiction.

Assume now that for some $c=e_1 \dots e_{2l}\in {\cal C}$ we have
$\psi(c) \in  < {\cal S}-\{\psi(c)\} >$. In this case there is a monomial
in one of the generators of $< {\cal S}-\{\psi(c)\} >$ which divides a monomial
in $\psi(c)$, i.e., there exists
$d=f_1 \dots f_{2r}\in {\cal C}$ such that after a cyclic permutation of
the edges of $d$ we have
$f_1=e_{2 i_1 -1}, f_3=e_{2 i_2 -1}, \dots, f_{2r-1}=e_{2 i_r -1}$.
But then if any of $f_{2}, f_4, \dots f_{2r}$ is not an edge in
$c$ then it must be a chord of $c$ and, therefore, all the edges of $d$ are in
$c$, implying that $c=d$, a contradiction.
\end{proof}

\begin{lem}\label{CI-Lemma}
Let $G$ be a graph with no triangles with the property that
any two cycles with no chords in $G$ share at most one edge.
\begin{enumerate}
\item If two cycles with no chords have a common edge then there
is no edge connecting them other than the common edge.
\item There exists an edge in at most one cycle with no chords.
\item If $e$ is an edge on precisely one cycle with no chords
then $e$ is not a chord of any cycle.
\item If $G$ is connected then there are precisely
$\nE(G)-\nV(G)+1$ cycles with no chords in $G$.
\end{enumerate}
\end{lem}

\begin{proof}
It is enough to prove the lemma for all the connected components
of $G$, so we will assume henceforth that $G$ is connected.
\begin{enumerate}

\item
Let $c_1,c_2\in \C(G)$ have the edge $(v_1,v_2)$ as a common edge.
Write
$$c_1=(v_1,v_2)(v_2,u_3)\dots (u_{r},v_1) \quad\quad(r>3)$$
and
$$c_2=(v_1,v_2)(v_2,w_3)\dots (w_{s},v_1) \quad\quad(s>3).$$

If there is an edge other than $(v_1,v_2)$ connecting
$c_1$ and $c_2$, since $c_1$ and $c_2$ have no chords we can pick
$3\le i \le r$ minimal such that there exists an edge
in $G$ connecting $u_i$ with a vertex of $c_2$, and we can pick
$3\le j\le s$ minimal such that $(u_i,w_j)$ is an edge.

We cannot have $i=j=3$ otherwise we would have a triangle in
$G$, and we may assume that $i>3$. The cycle
$$c_3=(v_2,u_3)(u_3,u_4)\dots(u_{i-1},u_i)(u_i,w_j)(w_j,w_{j-1})\dots
(w_3,v_2)$$
has no chords and $\#(c_3\cap c_1)>1$, a contradiction.

\item

Let $\mathfrak G$ be the bipartite graph whose vertices are
$$ \C(G) \cup
\left\{ e\in \E(G) | e \text{ is in some } c\in\C(G) \right\}$$
and whose edges are
$$ \left\{ (e,c) | e\in\E(G), c\in\C(G) \text{ and }
e \text{ is an edge of } c\right\} .$$
If any edge is in at least two cycles with no chords then
the degree of the vertices of $\mathfrak G$ is at least two and we can
pick a minimal cycle
$c_1,e_1,c_2,e_2,\dots,c_r,e_r$ in $\mathfrak G$, i.e.,
we produce
a sequence $(c_1,\dots,c_r)\subset \C(G)$ together with a sequence of
edges $e_1,\dots, e_r$ such that for all $1\le i < r$ we have
$e_i\in c_i\cap c_{i+1}$ and $e_r\in c_r \cap c_1$ and in addition
only consecutive (modulo $r$) cycles in this sequence have a common edge.

We first note that there is no edge connecting two vertices
in different $c_i,c_j$ other than one of $e_1,\dots,e_r$; if there were such
an edge $e$ then by part (1) of the lemma $i$ and $j$ are not consecutive
(modulo $r$.) After a cyclic permutation of the cycles we
may assume that $1\le i <j <r$ and write $e=(v_1,v_2)$ with
$v_1\in c_i$ and $v_2 \in c_j$. We can find a path $p$ between $v_1$ and
$v_2$ lying in
$(c_i \cup \dots \cup c_j) - \{e_i,e_{i+1},\dots, e_{j-1}\}$;
add to this path the edge $e$
to obtain a cycle $c$. If $c$ has chords replace it with another cycle
with no chords containing a sub-path of $p$ and an edge $e^\prime$ connecting
two vertices in $c_{i^\prime}$ and $c_{j^\prime}$ with
$i\le i^\prime < j^\prime \le j$. Thus we may assume that $c$ has no chords
and we may replace $c_1,\dots,c_r$ with a possibly shorter sequence
$c_i,c_{i+1},c$ implying that $r=3$. But when $r=3$ any two cycles are
consecutive and we are done by the first part of this lemma.

Consider the graph $H=(c_1 \cup \dots \cup c_r) - \{e_1,\dots,e_r\}$; $H$ has
at most
two connected components, one of which must be a cycle $c$
(one of these components may be a single vertex, but not both.)

Assume first that $H=c$. For any $e_i$ there is a path $p$ in $H$ connecting
the endpoints of $e_i$, and if we pick this path to have minimal length,
the cycle $c^\prime$ obtained by concatenating $p$ and $e_i$ has no chords.
But $p$ must have an edge in common with either $c_{i-1}$ or with $c_i$, and,
therefore, $c^\prime$ must share at least two edges with $c_{i-1}$ or with
$c_i$.

Consider now the case where $H$ has two connected components,
one of which is the cycle $c$.
We have shown that this cycle cannot have a chord, i.e.,
$c\in \C(G)$, and, therefore, $\#(c\cap c_i)\le 1$ for all
$1\le i\le r$. But since every edge in every $c_i$ except two are in
$H$ we must have  $\#(c\cap c_i)= 1$ for all
$1\le i\le r$. This immediately shows that both connected components
of $H$ are cycles and also that $r>3$ because $G$ has no
triangles.

Let the two connected components of $H$ be
$f_1,f_2, \dots, f_r$ and $g_1,g_2,\dots,g_r$ where $f_i,g_i\in c_i$
for all $1\le i \le r$ and consider the cycles
$$c^\prime=e_r, g_1, g_2, e_2, f_3, f_4, \dots, f_r  $$
and
$$c^{\prime\prime}=g_1,g_2,\dots, g_r .$$
These cycles have no chords and their intersection is $\{ g_1,g_2 \}$,
a contradiction.

\item Any chord is an edge of at least two cycles with no chords.

\item We proceed by induction on $\nE(G)-\nV(G)$.
If $\nE(G)-\nV(G)=-1$ then
$G$ is a tree and the claim is trivial.
Assume that $\nE(G)-\nV(G)\ge 0$ and pick an edge
 $e$ precisely  in one $c\in \C(G)$.
Consider the graph $H=G-\{e\}$; by removing the edge $e$ we removed from
$G$ one cycle with no chords, and since $e$ is not a chord of any cycle
in $G$, removing $e$ does not add any new cycles with no chords.
Thus $H$ has one less cycle with no chords than $G$ and
by the induction hypothesis $H$ has $\nE(G)-\nV(G)$ such cycles.

\end{enumerate}
\end{proof}

\begin{lem}\label{non-CI-Lemma}
Let $G$ be any graph.
\begin{enumerate}
\item $\#\C(G)\ge \nE(G)-\nV(G)+1$.
\item If $G$ has two cycles with no chords
with more than one common edge then
$\#\C(G)> \nE(G)-\nV(G)+1$.
\end{enumerate}
\end{lem}

\begin{proof}

\begin{enumerate}

\item
If $e$ is any edge in $G$ we denote by $G_{(e)}$ the graph obtained from
$G$ by ``shrinking'' $e$, i.e., by removing the edge $e$ and identifying its
endpoints. We also denote by $\Delta_e$ the number of triangles in
$G$ of which $e$ is an edge.

We have $\nE(G_{(e)})=\nE(G)-\Delta_e-1$ and $\nV(G_{(e)})=\nV(G)-1$.
We also have $\# \C(G_{(e)})=\# \C(G) - \Delta_e - \epsilon_e$
where $\epsilon_e\ge 0$ is the number of cycles with no chords in $G$ which
acquire a chord after shrinking $e$.

We can now use induction on $\nE(G)$:
\begin{gather*}
\#\C(G)=\#\C(G_{(e)})+\Delta_e + \epsilon_e \ge
\nE(G_{(e)})-\nV(G_{(e)})+1+\Delta_e = \\
\nE(G)-\Delta_e-1-\nV(G)+1+1+\Delta_e =
\nE(G)-\nV(G)+1 .
\end{gather*}

\item
Pick $c_1,c_2\in \C(G)$ such
that $\#(c_1\cap c_2)>1$. We can find a path
$$p=(u,w_1),(w_1,w_2),\dots,(w_l,v)$$
where $u,v$ are vertices in
$c_1$ and $w_1,\dots,w_l$ are vertices in $c_2-c_1$.
Note that $(u,v)$ cannot be an edge in $G$, otherwise, since
$c_2$ has no chords, $c_2$ would be the concatenation of $p$
and $(u,v)$ and would have only one edge in common with $c_1$.

We can now
shrink $G$ successively at all edges of $p$ but one. After this shrinking
$c_1$ will acquire a chord, thus at least one of the $\epsilon_e$'s
obtained in this process will be positive, and the inequality follows.

\end{enumerate}
\end{proof}

\begin{thm}\label{Main-theorem}
Let $G$ be a bipartite graph. $k[G]$ is a complete intersection
if and only if $G$ is a CI graph.
\end{thm}

\begin{proof}
If $G_1$ and $G_2$ are two disjoint graphs then
$k[G_1\cup G_2]=k[G_1] \otimes_k k[G_2]$, thus we may assume that
$G$ is connected.

$k[G]$ is a complete intersection
if and only if $K_G$ is generated by
$\nE(G)-\dim(k[G])=\nE(G)-\nV(G)+1$ elements
(cf. corollary 7.3.1 in \cite{V}) and
theorem \ref{CI-mingens} implies that
$k[G]$ is a complete intersection
if and only if  $\#\C(G)=\nE(G)-\nV(G)+1$; the result now follows
from lemmas \ref{CI-Lemma}(4) and \ref{non-CI-Lemma}.

\end{proof}

\begin{ex}\label{G_n-example}
Consider the graph $G_n$ with vertices
$\{ x,y,u_1,v_1,u_2,v_2,\dots,u_n,v_n\}$ and edges
$$\{(x,y)\} \cup
\{ (x,u_1),\dots,(x,u_n) \} \cup
\{ (y,v_1),\dots,(y,v_n) \} \cup
\{ (u_1,v_1), \dots, (u_n,v_n) \} .$$

$G_n$ is bipartite with $\#\C(G_n)=n$ and since
$\nE(G_n)-\nV(G_n)+1=3n+1-(2n+2)+1=n$ we conclude that $k[G_n]$ is a complete
intersection. Notice, however, that if $H_n$ is the graph obtained
from $G_n$ by removing the edge $(x,y)$
we have $\#\C(H_n)=\binom{n}{2}$ cycles with no chords, and, therefore,
$K_{H_n}$ is a prime ideal of height $n$ which is $\binom{n}{2}$--generated.
\end{ex}

\begin{thm}
Let $G$ be a graph as in lemma \ref{CI-Lemma}. Then $G$ is planar.
\end{thm}

\begin{proof}

The following proof is based on the proof of lemma 11.13(a) in \cite{H}.

Pick a counterexample $G$ with minimal $\nE(G)$; $G$ will necessarily be
a block and we may pick an edge
$e=(u_1,u_4)\in \E(G)$ lying in a unique $c\in \C(G)$.
We may shrink the edge $e$ in $G$ without affecting the hypothesis
of the theorem unless $c$ is a cycle of length four; we shall assume henceforth
that $c=(u_1,u_2),(u_2,u_3),(u_3,u_4),(u_4,u_1)$.

Let $H=G-\{e\}$; note that $H$ satisfies the hypothesis of the theorem and that
$u_1$ and $u_4$ must lie in different blocks
$B_1$ and $B_2$
of $H$ thus we may pick a cutpoint in all paths in $H$ from $u_1$ to
$u_4$ and with no loss of generality we may take this cutpoint to be
$u_2$.

Let $B_2^\prime=B_2 \cup \{(u_2,u_4),(u_2,u_3)\}$ and let
$B_2^{\prime\prime}=B_2 \cup \{(u_1,u_4),(u_1,u_2),(u_2,u_3)\}$
(note that the edge $(u_2,u_3)$ may have already been present in $B_2$.)
Clearly, $B_2^{\prime\prime}$ contains no triangles and since
the only cycle of $B_2^{\prime\prime}$ not in $B_2$ is $c$,
we see that $B_2^{\prime\prime}$ satisfies the hypothesis of the theorem.

If $B_2^{\prime\prime}\ne G$ we may deduce that it is planar, and
$B_2^\prime$, being homeomorphic to it, must also be planar. We
may then embed $H\cup \{(u_2,u_4)\}$ in the plane in such a way that
$(u_1,u_2)$ and $(u_2,u_4)$ are exterior edges; adding now the edge
$(u_1,u_4)$ will not affect the planarity of the graph, and we conclude that
$G$ is planar.

Assume now that $B_2^{\prime\prime}= G$. If $u_2$ and $u_4$ belong to
different blocks of $F=B_2\cup \{(u_2,u_3)\}$
then so do the edges $(u_2,u_3)$ and $(u_3,u_4)$ and we can
embed $F$ in the plane so that these edges bound the exterior face.
We can then add the edges $(u_1,u_2)$ and $(u_1,u_4)$ without affecting
the planarity. If $u_2$ and $u_4$ lie in the same block of $F$ we can find
minimal path $p$ in $F-\{u_3\}$ connecting $u_2$ with $u_4$. The cycle obtained
by concatenating $p$ with $(u_1,u_2)$ and $(u_1,u_4)$ has no chords
and is different from $c$, contradicting the fact that $(u_1,u_4)$ lies in
a unique cycle with no chords.
\end{proof}

\begin{cor}
Let $G$ be a connected CI graph. Then either
$G$ is a single edge or $\nE(G) \le 2(\nV(G)-2)$.
\end{cor}
\begin{proof}
Since $G$ must be planar and with no triangles, the result follows easily from
Euler's formula for planar graphs (see also
corollary 11.17(b) in \cite{H}.)
\end{proof}

\begin{rem}\

\begin{enumerate}

\item It is not hard to see that a bipartite outerplanar graph is a CI graph
but the reverse inclusion does not hold, e.g. the graph $G_n$ in example
\ref{G_n-example} is not outerplanar for $n\ge 3$ since it contains a subgraph
homeomorphic to $K_{2,3}$. Therefore the family of CI-graphs is strictly
between the families of
bipartite outerplanar graphs and
bipartite planar graphs.

\item When $G$ is not bipartite, $k[G]$ may be a complete intersection
without $G$ being planar. For example let $G$ be the following graph:
$$
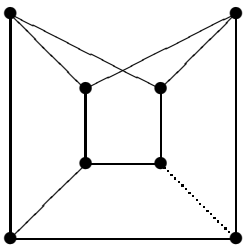
\begin{xy}
(10,20)*=0{\bullet}="1",
(20,20)*=0{\bullet}="2",
(30,0)*=0{\bullet}="3",
(0,30)*=0{\bullet}="4",
(30,30)*=0{\bullet}="5",
(0,0)*=0{\bullet}="6",
(10,10)*=0{\bullet}="7",
(20,10)*=0{\bullet}="8",
"1";"4"  **@@{-},
"1";"5"  **@@{-},
"1";"7"  **@@{-},
"2";"4"  **@@{-},
"2";"5"  **@@{-},
"2";"8"  **@@{-},
"8";"7"  **@@{-},
"6";"4"  **@@{-},
"6";"3"  **@@{-},
"3";"5"  **@@{-},
"6";"7"  **@@{-},
"3";"8"  **@@{.}
\end{xy}
$$
A computation with Macaulay2 (\cite{GS}) shows that $k[G]$ is a complete
intersection; the solid lines show a subgraph of $G$ homeomorphic to
$K_{3,3}$.
\end{enumerate}
\end{rem}

%%%%%%%%%%%%%%%%%%%%%%%%%%%%%%%%%%%%%%%%%%%%%%%%%%%%
\section{Algorithmic applications and some examples}
%%%%%%%%%%%%%%%%%%%%%%%%%%%%%%%%%%%%%%%%%%%%%%%%%%%%

In this section we will generalize theorem \ref{CI-mingens}
which will result in an algorithm for computing $\C(G)$.
Throughout this section we shall assume that $k[\E(G)]$ is equipped
with a monomial order so that for any closed walk $w$ of even length
$\psi(w)$ is well defined.

\begin{thm}
The elements of $\{\psi(c) | c\in \C(G) \text{ is an even cycle}\}$ form part of a minimal
set of generators for $K_G$.
\end{thm}

\begin{proof}
Let $W$ be a set of closed walks of even length such that
$\{ \psi(w) | w\in W \}$ is a minimal set of generators
for $K_G$ and let $c\in \C(G)$. We will
show that $c\in W$.

Since $\psi(c)\in K_G$ there exists a $w\in W$ and a monomial in $\psi(w)$
which divides $\psi^+(c)$.

If $w$ contains no odd cycles then the proof of theorem \ref{CI-mingens}
shows that $w\in \C(G)$ and that $w=c$.

If $w=(e_1,e_2,\dots,e_{2l})$ contains an odd cycle, say
$(e_1,e_2,\dots,e_{2r+1})$ then each of $\psi^+(w)$ and $\psi^-(w)$ is
divisible by one of $e_1 e_{2r+1}$ or $e_1 e_{2r+2}$. But this is impossible
since $\psi^+(c)$ is not divisible by any two edges sharing a common vertex.
\end{proof}

As a corollary we obtain an algorithm for producing $\C(G)$ as follows:
given a graph $G$ construct the ideal $I_G$ generated by
$$\{ e-u v | u,v\in V(G), e=(u,v)\in \E(G) \}
\subset R=k[\V(G),\E(G)] .$$
Using a lexicographic order in $R$ with $v > e$ for any
$v\in \V(G)$ and $e\in \E(G)$ compute a Gr\"obner basis for
$I_G$ and eliminate the variables corresponding to vertices of $G$.
The resulting set will contain a minimal subset of generators for
$K_G$; we can now pick those corresponding to $\C(G)$.

\begin{ex}
Let $G$ be the following graph:

$$
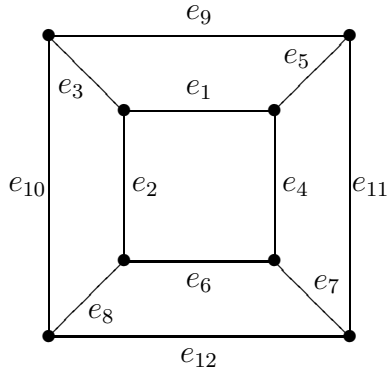
\begin{xy}
(10,30)*=0{\bullet}="1",
(30,30)*=0{\bullet}="2",
(30,10)*=0{\bullet}="3",
(10,10)*=0{\bullet}="4",
(0,40)*=0{\bullet}="5",
(40,40)*=0{\bullet}="6",
(40,0)*=0{\bullet}="7",
(0,0)*=0{\bullet}="8",
"1";"2"  **@@{-} ? *_!/8pt/{e_1},
"1";"4"  **@@{-} ? *_!/8pt/{e_2},
"1";"5"  **@@{-} ? *_!/8pt/{e_3},
"2";"3"  **@@{-} ? *_!/8pt/{e_4},
"2";"6"  **@@{-} ? *_!/8pt/{e_5},
"3";"4"  **@@{-} ? *_!/8pt/{e_6},
"3";"7"  **@@{-} ? *_!/8pt/{e_7},
"4";"8"  **@@{-} ? *_!/8pt/{e_8},
"5";"6"  **@@{-} ? *_!/8pt/{e_9},
"8";"5"  **@@{-} ? *_!/8pt/{e_{10}},
"6";"7"  **@@{-} ? *_!/8pt/{e_{11}},
"7";"8"  **@@{-} ? *_!/8pt/{e_{12}}
\end{xy}
$$

Applying the algorithm above using a lexicographical order in which
$e_1>e_2>\dots>e_{12}$ we obtain a Gr\"obner basis for
$K_G$ corresponding to the cycles:

$$
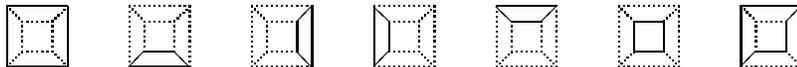
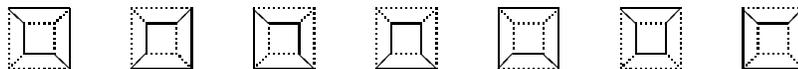
\begin{xy}
(2,6)*=0{}="1",
(6,6)*=0{}="2",
(6,2)*=0{}="3",
(2,2)*=0{}="4",
(0,8)*=0{}="5",
(8,8)*=0{}="6",
(8,0)*=0{}="7",
(0,0)*=0{}="8",
"1";"2"  **@@{.} ,
"1";"4"  **@@{.} ,
"1";"5"  **@@{.} ,
"2";"3"  **@@{.} ,
"2";"6"  **@@{.} ,
"3";"4"  **@@{.} ,
"3";"7"  **@@{.} ,
"4";"8"  **@@{.} ,
"5";"6"  **@@{-} ,
"8";"5"  **@@{-} ,
"6";"7"  **@@{-} ,
"7";"8"  **@@{-}
\end{xy}
\quad\quad
\begin{xy}
(2,6)*=0{}="1",
(6,6)*=0{}="2",
(6,2)*=0{}="3",
(2,2)*=0{}="4",
(0,8)*=0{}="5",
(8,8)*=0{}="6",
(8,0)*=0{}="7",
(0,0)*=0{}="8",
"1";"2"  **@@{.} ,
"1";"4"  **@@{.} ,
"1";"5"  **@@{.} ,
"2";"3"  **@@{.} ,
"2";"6"  **@@{.} ,
"3";"4"  **@@{-} ,
"3";"7"  **@@{-} ,
"4";"8"  **@@{-} ,
"5";"6"  **@@{.} ,
"8";"5"  **@@{.} ,
"6";"7"  **@@{.} ,
"7";"8"  **@@{-}
\end{xy}
\quad\quad
\begin{xy}
(2,6)*=0{}="1",
(6,6)*=0{}="2",
(6,2)*=0{}="3",
(2,2)*=0{}="4",
(0,8)*=0{}="5",
(8,8)*=0{}="6",
(8,0)*=0{}="7",
(0,0)*=0{}="8",
"1";"2"  **@@{.} ,
"1";"4"  **@@{.} ,
"1";"5"  **@@{.} ,
"2";"3"  **@@{-} ,
"2";"6"  **@@{-} ,
"3";"4"  **@@{.} ,
"3";"7"  **@@{-} ,
"4";"8"  **@@{.} ,
"5";"6"  **@@{.} ,
"8";"5"  **@@{.} ,
"6";"7"  **@@{-} ,
"7";"8"  **@@{.}
\end{xy}
\quad\quad
\begin{xy}
(2,6)*=0{}="1",
(6,6)*=0{}="2",
(6,2)*=0{}="3",
(2,2)*=0{}="4",
(0,8)*=0{}="5",
(8,8)*=0{}="6",
(8,0)*=0{}="7",
(0,0)*=0{}="8",
"1";"2"  **@@{.} ,
"1";"4"  **@@{-} ,
"1";"5"  **@@{-} ,
"2";"3"  **@@{.} ,
"2";"6"  **@@{.} ,
"3";"4"  **@@{.} ,
"3";"7"  **@@{.} ,
"4";"8"  **@@{-} ,
"5";"6"  **@@{.} ,
"8";"5"  **@@{-} ,
"6";"7"  **@@{.} ,
"7";"8"  **@@{.}
\end{xy}
\quad\quad
\begin{xy}
(2,6)*=0{}="1",
(6,6)*=0{}="2",
(6,2)*=0{}="3",
(2,2)*=0{}="4",
(0,8)*=0{}="5",
(8,8)*=0{}="6",
(8,0)*=0{}="7",
(0,0)*=0{}="8",
"1";"2"  **@@{-} ,
"1";"4"  **@@{.} ,
"1";"5"  **@@{-} ,
"2";"3"  **@@{.} ,
"2";"6"  **@@{-} ,
"3";"4"  **@@{.} ,
"3";"7"  **@@{.} ,
"4";"8"  **@@{.} ,
"5";"6"  **@@{-} ,
"8";"5"  **@@{.} ,
"6";"7"  **@@{.} ,
"7";"8"  **@@{.}
\end{xy}
\quad\quad
\begin{xy}
(2,6)*=0{}="1",
(6,6)*=0{}="2",
(6,2)*=0{}="3",
(2,2)*=0{}="4",
(0,8)*=0{}="5",
(8,8)*=0{}="6",
(8,0)*=0{}="7",
(0,0)*=0{}="8",
"1";"2"  **@@{-} ,
"1";"4"  **@@{-} ,
"1";"5"  **@@{.} ,
"2";"3"  **@@{-} ,
"2";"6"  **@@{.} ,
"3";"4"  **@@{-} ,
"3";"7"  **@@{.} ,
"4";"8"  **@@{.} ,
"5";"6"  **@@{.} ,
"8";"5"  **@@{.} ,
"6";"7"  **@@{.} ,
"7";"8"  **@@{.}
\end{xy}
\quad\quad
\begin{xy}
(2,6)*=0{}="1",
(6,6)*=0{}="2",
(6,2)*=0{}="3",
(2,2)*=0{}="4",
(0,8)*=0{}="5",
(8,8)*=0{}="6",
(8,0)*=0{}="7",
(0,0)*=0{}="8",
"1";"2"  **@@{.} ,
"1";"4"  **@@{.} ,
"1";"5"  **@@{.} ,
"2";"3"  **@@{-} ,
"2";"6"  **@@{-} ,
"3";"4"  **@@{-} ,
"3";"7"  **@@{.} ,
"4";"8"  **@@{-} ,
"5";"6"  **@@{-} ,
"8";"5"  **@@{-} ,
"6";"7"  **@@{.} ,
"7";"8"  **@@{.}
\end{xy}
$$

$$
\begin{xy}
(2,6)*=0{}="1",
(6,6)*=0{}="2",
(6,2)*=0{}="3",
(2,2)*=0{}="4",
(0,8)*=0{}="5",
(8,8)*=0{}="6",
(8,0)*=0{}="7",
(0,0)*=0{}="8",
"1";"2"  **@@{.} ,
"1";"4"  **@@{-} ,
"1";"5"  **@@{-} ,
"2";"3"  **@@{.} ,
"2";"6"  **@@{.} ,
"3";"4"  **@@{-} ,
"3";"7"  **@@{-} ,
"4";"8"  **@@{.} ,
"5";"6"  **@@{-} ,
"8";"5"  **@@{.} ,
"6";"7"  **@@{-} ,
"7";"8"  **@@{.}
\end{xy}
\quad\quad
\begin{xy}
(2,6)*=0{}="1",
(6,6)*=0{}="2",
(6,2)*=0{}="3",
(2,2)*=0{}="4",
(0,8)*=0{}="5",
(8,8)*=0{}="6",
(8,0)*=0{}="7",
(0,0)*=0{}="8",
"1";"2"  **@@{-} ,
"1";"4"  **@@{-} ,
"1";"5"  **@@{.} ,
"2";"3"  **@@{.} ,
"2";"6"  **@@{-} ,
"3";"4"  **@@{.} ,
"3";"7"  **@@{.} ,
"4";"8"  **@@{-} ,
"5";"6"  **@@{.} ,
"8";"5"  **@@{.} ,
"6";"7"  **@@{-} ,
"7";"8"  **@@{-}
\end{xy}
\quad\quad
\begin{xy}
(2,6)*=0{}="1",
(6,6)*=0{}="2",
(6,2)*=0{}="3",
(2,2)*=0{}="4",
(0,8)*=0{}="5",
(8,8)*=0{}="6",
(8,0)*=0{}="7",
(0,0)*=0{}="8",
"1";"2"  **@@{-} ,
"1";"4"  **@@{.} ,
"1";"5"  **@@{-} ,
"2";"3"  **@@{-} ,
"2";"6"  **@@{.} ,
"3";"4"  **@@{.} ,
"3";"7"  **@@{-} ,
"4";"8"  **@@{.} ,
"5";"6"  **@@{.} ,
"8";"5"  **@@{-} ,
"6";"7"  **@@{.} ,
"7";"8"  **@@{-}
\end{xy}
\quad\quad
\begin{xy}
(2,6)*=0{}="1",
(6,6)*=0{}="2",
(6,2)*=0{}="3",
(2,2)*=0{}="4",
(0,8)*=0{}="5",
(8,8)*=0{}="6",
(8,0)*=0{}="7",
(0,0)*=0{}="8",
"1";"2"  **@@{-} ,
"1";"4"  **@@{-} ,
"1";"5"  **@@{.} ,
"2";"3"  **@@{-} ,
"2";"6"  **@@{.} ,
"3";"4"  **@@{.} ,
"3";"7"  **@@{-} ,
"4";"8"  **@@{-} ,
"5";"6"  **@@{.} ,
"8";"5"  **@@{.} ,
"6";"7"  **@@{.} ,
"7";"8"  **@@{-}
\end{xy}
\quad\quad
\begin{xy}
(2,6)*=0{}="1",
(6,6)*=0{}="2",
(6,2)*=0{}="3",
(2,2)*=0{}="4",
(0,8)*=0{}="5",
(8,8)*=0{}="6",
(8,0)*=0{}="7",
(0,0)*=0{}="8",
"1";"2"  **@@{.} ,
"1";"4"  **@@{.} ,
"1";"5"  **@@{.} ,
"2";"3"  **@@{.} ,
"2";"6"  **@@{.} ,
"3";"4"  **@@{-} ,
"3";"7"  **@@{-} ,
"4";"8"  **@@{-} ,
"5";"6"  **@@{-} ,
"8";"5"  **@@{-} ,
"6";"7"  **@@{-} ,
"7";"8"  **@@{.}
\end{xy}
\quad\quad
\begin{xy}
(2,6)*=0{}="1",
(6,6)*=0{}="2",
(6,2)*=0{}="3",
(2,2)*=0{}="4",
(0,8)*=0{}="5",
(8,8)*=0{}="6",
(8,0)*=0{}="7",
(0,0)*=0{}="8",
"1";"2"  **@@{.} ,
"1";"4"  **@@{-} ,
"1";"5"  **@@{-} ,
"2";"3"  **@@{-} ,
"2";"6"  **@@{-} ,
"3";"4"  **@@{-} ,
"3";"7"  **@@{.} ,
"4";"8"  **@@{.} ,
"5";"6"  **@@{-} ,
"8";"5"  **@@{.} ,
"6";"7"  **@@{.} ,
"7";"8"  **@@{.}
\end{xy}
\quad\quad
\begin{xy}
(2,6)*=0{}="1",
(6,6)*=0{}="2",
(6,2)*=0{}="3",
(2,2)*=0{}="4",
(0,8)*=0{}="5",
(8,8)*=0{}="6",
(8,0)*=0{}="7",
(0,0)*=0{}="8",
"1";"2"  **@@{-} ,
"1";"4"  **@@{.} ,
"1";"5"  **@@{-} ,
"2";"3"  **@@{.} ,
"2";"6"  **@@{-} ,
"3";"4"  **@@{.} ,
"3";"7"  **@@{.} ,
"4";"8"  **@@{.} ,
"5";"6"  **@@{.} ,
"8";"5"  **@@{-} ,
"6";"7"  **@@{-} ,
"7";"8"  **@@{-}
\end{xy}
$$

The first ten elements give us $\C(G)$.

\end{ex}

\begin{ex}
The minimal generators of $K_G$ when $G$ is not bipartite can correspond
to quite complicated paths. Let $G$ be the following graph:

$$
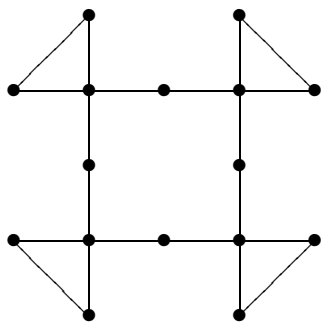
\begin{xy}
(30,30)*=0{\bullet}="1",
(30,20)*=0{\bullet}="2",
(30,10)*=0{\bullet}="3",
(20,10)*=0{\bullet}="4",
(10,10)*=0{\bullet}="5",
(10,20)*=0{\bullet}="6",
(10,30)*=0{\bullet}="7",
(20,30)*=0{\bullet}="8",
(30,40)*=0{\bullet}="9",
(40,30)*=0{\bullet}="10",
(40,10)*=0{\bullet}="11",
(30,0)*=0{\bullet}="12",
(10,0)*=0{\bullet}="13",
(0,10)*=0{\bullet}="14",
(0,30)*=0{\bullet}="15",
(10,40)*=0{\bullet}="16",
"1";"2"  **@@{-} ,
"2";"3"  **@@{-} ,
"3";"4"  **@@{-} ,
"4";"5"  **@@{-} ,
"5";"6"  **@@{-} ,
"6";"7"  **@@{-} ,
"7";"8"  **@@{-} ,
"8";"1"  **@@{-} ,
"1";"9"  **@@{-} ,
"1";"10"  **@@{-} ,
"9";"10"  **@@{-} ,
"3";"11"  **@@{-} ,
"3";"12"  **@@{-} ,
"11";"12"  **@@{-} ,
"5";"13"  **@@{-} ,
"5";"14"  **@@{-} ,
"13";"14"  **@@{-} ,
"7";"15"  **@@{-} ,
"7";"16"  **@@{-} ,
"15";"16"  **@@{-}
\end{xy}
$$

There are twenty minimal generators of $K_G$
corresponding  to the following walks (up to symmetry):

$$
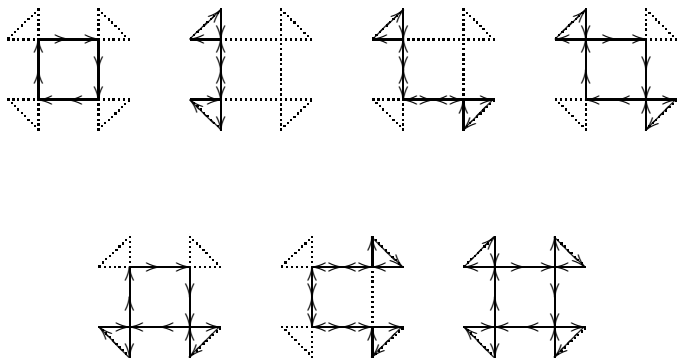
\begin{xy}
(12,12)*=0{}="1",(12,8)*=0{}="2",(12,4)*=0{}="3",
(8,4)*=0{}="4",(4,4)*=0{}="5",(4,8)*=0{}="6",
(4,12)*=0{}="7",(8,12)*=0{}="8",(12,16)*=0{}="9",
(16,12)*=0{}="10",(16,4)*=0{}="11",(12,0)*=0{}="12",
(4,0)*=0{}="13",(0,4)*=0{}="14",(0,12)*=0{}="15",
(4,16)*=0{}="16",
"1";"2"  **@@{.} ,"2";"3"  **@@{.} ,"3";"4"  **@@{.} ,"4";"5"  **@@{.} ,
"5";"6"  **@@{.} ,"6";"7"  **@@{.} ,"7";"8"  **@@{.} ,"8";"1"  **@@{.} ,
"1";"9"  **@@{.} ,"1";"10"  **@@{.} ,"9";"10"  **@@{.} ,"3";"11"  **@@{.} ,
"3";"12"  **@@{.} ,"11";"12"  **@@{.} ,"5";"13"  **@@{.} ,"5";"14"  **@@{.} ,
"13";"14"  **@@{.} ,"7";"15"  **@@{.} ,"7";"16"  **@@{.} ,"15";"16"  **@@{.},
@+"1"\PATH ~={**\dir{-}?>*\dir{>}}
 '"2" '"3" '"4" '"5" '"6" '"7" '"8" "1"
\end{xy}
\quad\quad
\begin{xy}
(12,12)*=0{}="1",(12,8)*=0{}="2",(12,4)*=0{}="3",
(8,4)*=0{}="4",(4,4)*=0{}="5",(4,8)*=0{}="6",
(4,12)*=0{}="7",(8,12)*=0{}="8",(12,16)*=0{}="9",
(16,12)*=0{}="10",(16,4)*=0{}="11",(12,0)*=0{}="12",
(4,0)*=0{}="13",(0,4)*=0{}="14",(0,12)*=0{}="15",
(4,16)*=0{}="16",
"1";"2"  **@@{.} ,"2";"3"  **@@{.} ,"3";"4"  **@@{.} ,"4";"5"  **@@{.} ,
"5";"6"  **@@{.} ,"6";"7"  **@@{.} ,"7";"8"  **@@{.} ,"8";"1"  **@@{.} ,
"1";"9"  **@@{.} ,"1";"10"  **@@{.} ,"9";"10"  **@@{.} ,"3";"11"  **@@{.} ,
"3";"12"  **@@{.} ,"11";"12"  **@@{.} ,"5";"13"  **@@{.} ,"5";"14"  **@@{.} ,
"13";"14"  **@@{.} ,"7";"15"  **@@{.} ,"7";"16"  **@@{.} ,"15";"16"  **@@{.}
@+"5"\PATH ~={**\dir{-}?>*\dir{>}}
 '"13" '"14" '"5" '"6" '"7" '"15" '"16" '"7" '"6" "5"
\end{xy}
\quad\quad
\begin{xy}
(12,12)*=0{}="1",(12,8)*=0{}="2",(12,4)*=0{}="3",
(8,4)*=0{}="4",(4,4)*=0{}="5",(4,8)*=0{}="6",
(4,12)*=0{}="7",(8,12)*=0{}="8",(12,16)*=0{}="9",
(16,12)*=0{}="10",(16,4)*=0{}="11",(12,0)*=0{}="12",
(4,0)*=0{}="13",(0,4)*=0{}="14",(0,12)*=0{}="15",
(4,16)*=0{}="16",
"1";"2"  **@@{.} ,"2";"3"  **@@{.} ,"3";"4"  **@@{.} ,"4";"5"  **@@{.} ,
"5";"6"  **@@{.} ,"6";"7"  **@@{.} ,"7";"8"  **@@{.} ,"8";"1"  **@@{.} ,
"1";"9"  **@@{.} ,"1";"10"  **@@{.} ,"9";"10"  **@@{.} ,"3";"11"  **@@{.} ,
"3";"12"  **@@{.} ,"11";"12"  **@@{.} ,"5";"13"  **@@{.} ,"5";"14"  **@@{.} ,
"13";"14"  **@@{.} ,"7";"15"  **@@{.} ,"7";"16"  **@@{.} ,"15";"16"  **@@{.}
@+"3"\PATH ~={**\dir{-}?>*\dir{>}}
 '"11" '"12" '"3" '"4" '"5" '"6" '"7" '"15" '"16" '"7" '"6" '"5" '"4" "3"
\end{xy}
\quad\quad
\begin{xy}
(12,12)*=0{}="1",(12,8)*=0{}="2",(12,4)*=0{}="3",
(8,4)*=0{}="4",(4,4)*=0{}="5",(4,8)*=0{}="6",
(4,12)*=0{}="7",(8,12)*=0{}="8",(12,16)*=0{}="9",
(16,12)*=0{}="10",(16,4)*=0{}="11",(12,0)*=0{}="12",
(4,0)*=0{}="13",(0,4)*=0{}="14",(0,12)*=0{}="15",
(4,16)*=0{}="16",
"1";"2"  **@@{.} ,"2";"3"  **@@{.} ,"3";"4"  **@@{.} ,"4";"5"  **@@{.} ,
"5";"6"  **@@{.} ,"6";"7"  **@@{.} ,"7";"8"  **@@{.} ,"8";"1"  **@@{.} ,
"1";"9"  **@@{.} ,"1";"10"  **@@{.} ,"9";"10"  **@@{.} ,"3";"11"  **@@{.} ,
"3";"12"  **@@{.} ,"11";"12"  **@@{.} ,"5";"13"  **@@{.} ,"5";"14"  **@@{.} ,
"13";"14"  **@@{.} ,"7";"15"  **@@{.} ,"7";"16"  **@@{.} ,"15";"16"  **@@{.}
@+"3"\PATH ~={**\dir{-}?>*\dir{>}}
 '"11" '"12" '"3" '"4" '"5" '"6" '"7" '"15" '"16" '"7" '"8" '"1" '"2" "3"
\end{xy}
$$

$$
\begin{xy}
(12,12)*=0{}="1",(12,8)*=0{}="2",(12,4)*=0{}="3",
(8,4)*=0{}="4",(4,4)*=0{}="5",(4,8)*=0{}="6",
(4,12)*=0{}="7",(8,12)*=0{}="8",(12,16)*=0{}="9",
(16,12)*=0{}="10",(16,4)*=0{}="11",(12,0)*=0{}="12",
(4,0)*=0{}="13",(0,4)*=0{}="14",(0,12)*=0{}="15",
(4,16)*=0{}="16",
"1";"2"  **@@{.} ,"2";"3"  **@@{.} ,"3";"4"  **@@{.} ,"4";"5"  **@@{.} ,
"5";"6"  **@@{.} ,"6";"7"  **@@{.} ,"7";"8"  **@@{.} ,"8";"1"  **@@{.} ,
"1";"9"  **@@{.} ,"1";"10"  **@@{.} ,"9";"10"  **@@{.} ,"3";"11"  **@@{.} ,
"3";"12"  **@@{.} ,"11";"12"  **@@{.} ,"5";"13"  **@@{.} ,"5";"14"  **@@{.} ,
"13";"14"  **@@{.} ,"7";"15"  **@@{.} ,"7";"16"  **@@{.} ,"15";"16"  **@@{.}
@+"3"\PATH ~={**\dir{-}?>*\dir{>}}
 '"11" '"12" '"3" '"4" '"5" '"13" '"14" '"5" '"6" '"7" '"8" '"1" '"2" "3"
\end{xy}
\quad\quad
\begin{xy}
(12,12)*=0{}="1",(12,8)*=0{}="2",(12,4)*=0{}="3",
(8,4)*=0{}="4",(4,4)*=0{}="5",(4,8)*=0{}="6",
(4,12)*=0{}="7",(8,12)*=0{}="8",(12,16)*=0{}="9",
(16,12)*=0{}="10",(16,4)*=0{}="11",(12,0)*=0{}="12",
(4,0)*=0{}="13",(0,4)*=0{}="14",(0,12)*=0{}="15",
(4,16)*=0{}="16",
"1";"2"  **@@{.} ,"2";"3"  **@@{.} ,"3";"4"  **@@{.} ,"4";"5"  **@@{.} ,
"5";"6"  **@@{.} ,"6";"7"  **@@{.} ,"7";"8"  **@@{.} ,"8";"1"  **@@{.} ,
"1";"9"  **@@{.} ,"1";"10"  **@@{.} ,"9";"10"  **@@{.} ,"3";"11"  **@@{.} ,
"3";"12"  **@@{.} ,"11";"12"  **@@{.} ,"5";"13"  **@@{.} ,"5";"14"  **@@{.} ,
"13";"14"  **@@{.} ,"7";"15"  **@@{.} ,"7";"16"  **@@{.} ,"15";"16"  **@@{.}
@+"1"\PATH ~={**\dir{-}?>*\dir{>}}
'"9" '"10" '"1" '"8" '"7" '"6" '"5" '"4" '"3" '"11" '"12" '"3"
'"4" '"5" '"6" '"7" '"8" "1"
\end{xy}
\quad\quad
\begin{xy}
(12,12)*=0{}="1",(12,8)*=0{}="2",(12,4)*=0{}="3",
(8,4)*=0{}="4",(4,4)*=0{}="5",(4,8)*=0{}="6",
(4,12)*=0{}="7",(8,12)*=0{}="8",(12,16)*=0{}="9",
(16,12)*=0{}="10",(16,4)*=0{}="11",(12,0)*=0{}="12",
(4,0)*=0{}="13",(0,4)*=0{}="14",(0,12)*=0{}="15",
(4,16)*=0{}="16",
"1";"2"  **@@{.} ,"2";"3"  **@@{.} ,"3";"4"  **@@{.} ,"4";"5"  **@@{.} ,
"5";"6"  **@@{.} ,"6";"7"  **@@{.} ,"7";"8"  **@@{.} ,"8";"1"  **@@{.} ,
"1";"9"  **@@{.} ,"1";"10"  **@@{.} ,"9";"10"  **@@{.} ,"3";"11"  **@@{.} ,
"3";"12"  **@@{.} ,"11";"12"  **@@{.} ,"5";"13"  **@@{.} ,"5";"14"  **@@{.} ,
"13";"14"  **@@{.} ,"7";"15"  **@@{.} ,"7";"16"  **@@{.} ,"15";"16"  **@@{.}
@+"1"\PATH ~={**\dir{-}?>*\dir{>}}
'"9" '"10" '"1" '"2" '"3" '"11" '"12" '"3" '"4" '"5" '"13" '"14"
'"5" '"6" '"7" '"15" '"16" '"7" '"8" "1"
\end{xy}
$$

Notice that the last generator corresponds to the Euler path in $G$.
It is possible to generalize this example to obtain Eulerian graphs
in which the Euler paths correspond to minimal generators of $K_G$  and
where there are minimal generators of $K_G$ corresponding to closed walks
containing an arbitrarily large number of odd cycles.

\end{ex}

\section*{Acknowledgment}

I would like to express my gratitude to Victor Neumann-Lara for our pleasant
discussions on
Graph Theory and his many useful suggestions.

%%%%%%%%%%%%%%%%%%%%%%%%%%%%%%%%%%%%%%%%%%%%%%%%%%%%%%%%%%%%%%%%%%%%%%%%%%%%%%%

\end{document}